\documentclass[a4paper,12pt]{article}

\newtheorem{theorem}{Theorem}
\newtheorem{cor}[theorem]{Corollary}
\newtheorem{lemma}[theorem]{Lemma}

\begin{document}

\begin{center}
{\large\bf On decomposition of commutative Moufang groupoids}
\end{center}

\begin{center}
{\it B.\,V.\,Novikov}
\end{center}

It is well-known that the multiplicative groupoid of an alternative/Jordan algebra
satisfies Moufang identities \cite{bel, jac}. Therefore it seems interesting to study the
structure of such groupoids. In this note we apply to Moufang groupoids an approach which
is widespread in Semigroup Theory  --- decomposition into a semilattice of subsemigroups
\cite{c-p}.

\bigskip

We shall call a groupoid with the identity
\begin{equation}\label{eq-1}
(xy)(zx)=(x(yz))x
\end{equation}
by {\it a Moufang groupoid}. {\bf Everywhere in this article $M$
denotes a commutative Moufang groupoid.}

\begin{theorem}\label{th-1}
If $M$ consists of idempotents, then it is an semilattice.
\end{theorem}
{\bf Proof.} Under assumption of the theorem it follows from
(\ref{eq-1}) for $y=z$
\begin{equation}\label{eq-2}
(xy)x=xy
\end{equation}
Applying (\ref{eq-2}) to the right part of (\ref{eq-1}), we get:
\begin{equation}\label{eq-3}
x(yz)=(xy)(xz)
\end{equation}

Define a binary relation $\le$ on $M$ :
$$
a\le b \Longleftrightarrow ab=a
$$
and show that it is a partial order.

Indeed, the reflexivity follows from idempotentness, the symmetry
follows from commutativity. Let $a\le b\le c$. Then
$$
ac=(ab)c=(ac)(bc)=(ac)b=(ab)(bc)=ab=a,
$$
i.\,e. $a\le c$.

Further, $ab$ is a greatest lower bound for the pair $\{a,b\}$.
Really, $ab\le a$, $ab\le b$ by (\ref{eq-2}). Suppose that $x\le
a$, $x\le b$. Then $(ab)x=(ax)(bx)=x\cdot x=x$, i.\,e. $x\le ab$.
\rule{7pt}{7pt}

\begin{lemma}\label{th-2}
$M$ is a groupoid with associative powers.
\end{lemma}
{\bf Proof.} For $a\in M$ we shall denote by $a^{(n)}$ an arbitrary term of the length
$n\ge 1$, all letters of which are $a$. If all such terms coincide in $M$, we denote them
by $a^n$.

We use the induction on length of the term. Let $a^{(k)}=a^k$ for
any $k<n$ (for $k=3$ this follows from commutativity). Consider
some term $a^{(n)}$. It can be written in the form
$a^{(n)}=a^{(k)}a^{(l)}$, where $k,l\ge 1$ and $k+l=n$; in view
of commutativity one can assume that $k\le l$.

Suppose that $k\ge 2$. Then under hypothesis of the induction
$$
a^{(n)}=a^{(k)}a^{(l)}=a^{k}a^{l}=(aa^{k-1})(aa^{l-1})=(a(a^{k-1}a^{l-1}))a=
(aa^{n-2})a=aa^{n-1}.
$$
Hence all terms of the form $a^{(n)}$ are equal. \rule{7pt}{7pt}

We denote by $L_a$ the left translation corresponding to an
element $a$: $L_ab=ab$. From (\ref{eq-1}) we have:
$$
(xy)^2=L_x^2y^2.
$$

We generalize this identity:
\begin{lemma}\label{th-3}
$(ab)^{2^n}=L_a^{2^n} b^{2^n}$ for any $a,b\in M,$ $n\ge 0$ (here powers are defined
correctly in view of Lemma \ref{th-2}).
\end{lemma}
{\bf Proof} by induction on $n$. Assume that for $n$ the
statement is faithful and prove it for $n+1$:
$$
(ab)^{2^{n+1}}=[(ab)^2]^{2^n}=[a(ab^2)]^{2^n}=L_a^{2^n}
(ab^2)^{2^n}=L_a^{2^n}L_a^{2^n} b^{2^{n+1}}=L_a^{2^{n+1}}
b^{2^{n+1}}. \ \rule{7pt}{7pt}
$$
\begin{cor}\label{th-4}
$(L_{a_1}\ldots L_{a_{k-1}}a_k)^{2^n}=L_{a_1}^{2^n}\ldots
L_{a_{k-1}}^{2^n}a_k^{2^n}$. \rule{7pt}{7pt}
\end{cor}

Further we shall need one more equality for translations:

\begin{lemma}\label{th-5}
$L_a^{2n}L_b=L_{L_a^nb}L_a^n$ for any $a,b\in M,$ $n\ge 1$.
\end{lemma}
{\bf Proof.} For $n=1$ this statement coincides with
(\ref{eq-1}). General case is obtained by induction on $n$.
\rule{7pt}{7pt}

Let $I_a$ be denoted the principal ideal, generated by $a\in M$.
It is clear that each element from $I_a$ can be written in the
form $L_{x_1}\ldots L_{x_{k-1}}L_{x_k}a$.

Define relations $\rho$ and $\sigma$:
\begin{equation}\label{eq-4}
a\rho b \Longleftrightarrow \exists n\ge 1\quad a^n\in I_b,
\end{equation}
\begin{equation}\label{eq-5}
a\sigma b \Longleftrightarrow a\rho b\quad  \& \quad b\rho a.
\end{equation}

\begin{lemma}\label{th-6}
$\sigma$ is a congruence.
\end{lemma}
{\bf Proof.} Reflexivity and symmetry are obvious, it is enough to
check transitivity and stability of $\rho$. Note that one can
assume in the definition of $\rho$ that $n$ is the power of the
two.

Let $a\rho b$, $b\rho c$, i.\,e.
$$
a^{2^m}=L_{x_1}\ldots L_{x_k}b,\quad b^{2^n}=L_{y_1}\ldots
L_{y_l}c.
$$
By Corollary \ref{th-4}
$$
a^{2^{m+n}}=L_{x_1}^{2^n}\ldots L_{x_k}^{2^n}b^{2^n}
=L_{x_1}^{2^n}\ldots L_{x_k}^{2^n}L_{y_1}\ldots L_{y_l}c \in I_c,
$$
so $\rho$ is transitive.

Now let $a\rho b$, i.\,e. $a^{2^n}=L_{x_1}\ldots L_{x_k}b$, and
$c\in M$.

1) Suppose that $k\le n$. Then using several times Lemma
\ref{th-5}, we get for some $u_1,\ldots,u_k\in M$:
$$
(ca)^{2^n}=L_c^{2^n}a^{2^n}=L_c^{2^n}L_{x_1}\ldots L_{x_k}b =
L_{u_1}\ldots L_{u_k}L_c^{2^{n-k}}b \in I_{cb}.
$$

2) Let $k>n$. Then $a^{2^{n+k+1}}=L_yL_{x_1}\ldots L_{x_k}b$,
where $y=a^{2^{n+k+1}-2^n}$. Since $k+1<n+k+1$, we get the case
1). Consequently, $ca\rho cb$. \rule{7pt}{7pt}

\begin{lemma}\label{th-7}
$M/\sigma$ is a semilattice.
\end{lemma}
{\bf Proof.} Obviously, $a\sigma a^2$ for any $a\in M$. So
$M/\sigma$ is an idempotent groupoid. By Theorem \ref{th-1} it is
a semilattice. \rule{7pt}{7pt}

Now let us to consider the structure of $\sigma$-classes (of
course, they are subgroupoids).

Like to Theory of Semigroups, we call a groupoid $M$ {\it
Archimedean} if $a\sigma b$ for any $a,b\in M$, where $\sigma$ is
defined by the conditions (\ref{eq-4}) and (\ref{eq-5}). It is
clear that an Archime\-dean groupoid is indecomposable into a
semilattice of subgroupoids.

\begin{lemma}\label{th-8}
Let $\sigma$ be a congruence on $M$, defined by conditions
(\ref{eq-4}) and (\ref{eq-5}). Then each $\sigma$-class is
Archimedean.
\end{lemma}
{\bf Proof.} Let $N$ is a $\sigma$-class, $a,b\in N$. Then
\begin{equation}\label{eq-6}
a^n=L_{x_1}\ldots L_{x_k}b
\end{equation}
for some $n>0$, $x_1,\ldots,x_k\in M$. We need to prove that in
the equality (\ref{eq-6}) elements $x_1,\ldots,x_k$ can be chosen
from $N$.

From (\ref{eq-6}) and Lemma \ref{th-5} we have:
$$
a^{n+2^k}=L_a^{2^k}L_{x_1}\ldots L_{x_k}b= L_{L_a^{2^{k-1}}x_1}
L_{L_a^{2^{k-2}}x_2}\ldots L_{L_ax_k}b.
$$

Show that for any $i\le k$ the element $y_i=L_a^{2^{k-i}}x_i$ is
contained in $N$. Indeed, since $y_i=a(L_a^{2^{k-i}-1}x_i)$, then
$y_i\rho a$. On the other hand,
$$
a^{n+2^k}=L_{y_1}\ldots L_{y_k}b=L_{y_1}\ldots
L_{y_{i-1}}[(L_{y_{i+1}}\ldots L_{y_k}b)y_i],
$$
whence $a\rho y_i$. Thereby, $a\sigma y_i$, i.\,e. $y_i\in N$.
\rule{7pt}{7pt}

The final result:

\begin{theorem}\label{th-9}
A commutative Moufang groupoid is a semilattice of Archime\-dean
groupoids. \rule{7pt}{7pt}
\end{theorem}

\bigskip

\noindent{\bf Example.} Let a finite semigroup $S$ yield the identity $ab=a$ (a left zero
semigroup), $F$ be a field, ${\rm char}\,F\ne 2$, $A=FS$ be the semigroup algebra. $A$ is
is a Jordan algebra with respect to the operation $x*y=\frac{1}{2}(xy+yx)$. Denote by
$A^*$ its multiplicative groupoid (as is well-known it is Moufang and commutative
\cite{jac}).

The operation in $A^*$ can be written as follows. For $x=\sum_{a\in S}\alpha_a a\in A^*$,
$\alpha_a \in F$, denote $|x|=\sum_{a\in S}\alpha_a$. Then
$$
x*y=\frac{1}{2}(|x|y+|y|x)
$$

From here $ x^{*n}=|x|^{n-1}x$ and $|xy|=|x||y|$. In particular, $x\in {\rm Rad}\,A^*$
iff $|x|=0$.

Evidently, all elements from ${\rm Rad}\,A^*$ constitute one $\sigma$-class. On the other
hand, if $x,y\not\in {\rm Rad}\,A^*$ then they divide one another. To make sure that, it
is enough to put
$$
t=\frac{1}{|x|^2}(2|x|y-|y|x);
$$
then $y=x*t$. Thus $A^*={\rm Rad}\,A^*\cup (A^*\setminus {\rm Rad}\,A^*)$ is the
decomposition of $A^*$ into Archimedean components.

\bigskip

Finally we discuss some problems which arise here.

1. For loops the identity \ref{eq-1} (central Moufang identity) is equal to each of ones
$x(y(xz))=((xy)x)z$ and  $((zx)y)x=z(x(yx))$ (left and right Moufang identities). This is
valid for multiplicative groupoids of Jordan algebras as well, but not in the general
case. So we can consider left and right Moufang (commutative) groupoids. Are there
similar decompositions for them?

2. Is there Archimedean decomposition in noncommutative situation? This is the case for
semigroups \cite{b-c}.

3. What can one say about the structure of an Archimedean component? For instance, can it
contain more than one idempotent (cf. \cite{c-p}, Ex.4.3.2)?

\end{document}